\documentclass{article}

\usepackage{amsthm}
\usepackage{amsmath}
\usepackage{amssymb}

\usepackage{graphicx}
\usepackage{framed}
\usepackage{subfigure}
\usepackage{float}

\usepackage{hyperref}

\usepackage{tikz}

\theoremstyle{definition}

\newtheorem{thm}{Theorem}[section]

\newtheorem{defn}[thm]{Definition}
\newtheorem{rem}[thm]{Remark}

\newcommand{\QED}{\ensuremath{\null\hfill\square}}

\newcommand{\R}{\mathbb{R}}
\newcommand{\Rplus}{\mathbb{R}_{\geq 0}}

\newcommand{\dee}{\mathrm{d}}

\newcommand{\SD}{\mathrm{SD}}
\newcommand{\SE}{\mathrm{SE}}

\newcommand{\var}{\mathrm{Var}}
\newcommand{\covar}{\mathrm{Cov}}
\newcommand{\mean}{\mathrm{mean}}

\newcommand{\bigO}{\mathrm{O}}

\usepackage{xcolor}

\usepackage{hyperref}
\title{Efficiently Approximating Spread Dimension with High Confidence}

\author{Kevin Dunne\footnote{\url{kevin.dunne[at]mailbox.org}}}

\begin{document}
\maketitle

\begin{abstract}
The concepts of \emph{spread} and \emph{spread dimension} of a metric 
space were introduced by Willerton in the context of quantifying 
biodiversity of ecosystems. In previous work, we developed the 
theoretical basis for applications of spread dimension as an 
intrinsic dimension estimator. In this paper we introduce the 
\emph{pseudo spread dimension} which is an efficient approximation of 
spread dimension, and we derive a formula for the standard error associated with this approximation.
\end{abstract}

\section{Introduction}\label{sec:Introduction}

Willerton introduced the \emph{spread} of a metric space as a measure of the ``size'' of a metric space in the context of quantifying biodiversity of 
ecosystems \cite{Willerton2013:Spread} and it is closely related to the \emph{magnitude} of a metric space \cite{Leinster2013:Magnitude,LeinsterMeckes2017:Magnitude}. Willerton showed that spread characterises a notion of dimension based on 
the growth rate of the spread as a space is scaled -- the \emph{spread dimension} -- and 
provided empirical evidence of a relationship between spread dimension
and Minkowski dimension \cite[Section 4]{Willerton2013:Spread}.

In \cite{Dunne2023:SpreadDimension} we built
upon the results of Willerton, and further developed the concept of spread dimension, particularly in the 
context of Riemannian manifolds and finite subsets of Riemannian 
manifolds. We framed these results as providing the theoretical basis of spread dimension as an intrinsic dimension estimator and demonstrated the practical applications with real and synthetic data \cite[Section 4]{Dunne2023:SpreadDimension}.

In the present work, we introduce the \emph{$S$-pseudo spread dimension} of a metric 
space for an arbitrary subset $S \subset X$ which can be computed in space and time 
complexity $\bigO(|S|\times|X|)$, compared to $\bigO(|X|^2)$ for the spread dimension. We show the pseudo spread 
dimension approximates the spread dimension, and we derive a formula for its standard error and confidence intervals.

Algorithms computing pseudo spread dimension and corresponding confidence intervals has been implemented in Python\footnote{\url{https://github.com/dk-gh/spread_dimension}}.

\subsection*{Defining Spread Dimension}

Spread is defined for any compact metric space $(X,d)$ equipped with probability measure. The spread of a metric space is a strictly positive real value, and since a metric 
space can be scaled by any constant factor $t \in \R_{>0}$, to define a new metric space $td(x,y)$, the spread therefore yields 
a one-parameter family of values associated with the underlying space.

\begin{defn}
Let $(X,d)$ be a compact metric space equipped with probability measure $\mu$.
The \emph{spread} $\sigma_d(t): \Rplus \rightarrow \Rplus$ is the
function defined by
\begin{equation}\label{eq:definition_of_spread}
\sigma_d(t) = \int_{x\in X} \frac{\dee \mu(x)}{\int_{y\in X}e^{-t d(x,y)}\dee\mu(y)}
\end{equation}
For a given value of $t$, we call the value $\sigma_d(t)$ the \emph{spread of} $(X,d)$ \emph{at scale} $t$.
\end{defn}

Spread is defined for a broad class of metric spaces, including all 
finite metric spaces and compact Riemannian manifolds. If $(X,d)$ is finite and equipped with the uniform probability distribution, then the spread is given by
\begin{equation}\label{eq:FiniteSpread}
\sigma_{d}(t) = \sum\limits_{x \in X}\frac{1}{\sum\limits_{y \in X}e^{-td(x,y)}}
\end{equation}

For a finite metric space $(X,d)$ the value $\sigma_d(t)$ can be interpreted as the number of ``effective'' points the space has at a given scale $t$. At scales close to zero the space resembles a single point with $\sigma_d(t)\rightarrow 1$ as $t\rightarrow 0$, while at large scales the space resembles $|X|$ discrete points with $\sigma_d(t) \rightarrow |X|$ as $t \rightarrow \infty$. One of the key 
insights of Willerton is that the 
growth of the spread for intermediate values of $t$ encodes 
geometric information about the space; it characterises a notion of 
dimension \cite[Section 4]{Willerton2013:Spread}.

In \cite{Willerton2013:Spread} Willerton defines the \emph{instantaneous 
growth}
$\mathcal{G}f(t)$ of a function $f(t)$ where

\begin{equation}\label{eq:definitionGrowthDimension}
\mathcal{G}f(t) = \frac{\dee \ln (f(t))}{\dee \ln(t)} =
\frac{t}{f(t)}\frac{\dee}{\dee t}f(t) \end{equation} with the second
expression following from an application of the chain rule.

\begin{defn} For $(X,d)$ a compact metric space with probability measure $\mu$, the \emph{instantaneous spread dimension},
or simply the \emph{spread dimension} is defined as the instantaneous growth of the spread dimension $\mathcal{G}\sigma_d(t)$. For each $t\in \Rplus$ we call the value $\mathcal{G}\sigma_d(t)$ the \emph{spread dimension of} $(X,d)$ \emph{at scale} $t$.
\end{defn}

Intuitively, for finite metric spaces, the spread dimension at scale $t$ can be interpreted as the dimension the space ``resembles'' at that scale: at small scales with $t \rightarrow 0$ the space resembles a single point, with dimension zero; and at large scales $t \rightarrow \infty$ the space resembles $|X|$ discrete points, also with dimension zero. If a sample is large enough, then for intermediate values of $t$ the spread dimension typically peaks or plateaus around the value corresponding with the \emph{intrinsic dimension} of the space, where intrinsic dimension is understood in the usual sense from data science and machine learning \cite{Camastra2016:IntrinsicDimension,CampadelliEtAl2015:IntrinsicDimension}.

\subsection*{Spread Dimension as an Intrinsic Dimension Estimator}

The relationship between spread dimension and the existing notion of intrinsic dimension can be made rigorous in the case of finite subsets of Riemannian manifolds, and in particular, submanifolds of Euclidean space $M\subset \R^k$ equipped with the Euclidean distance function -- i.e. data that satisfies the \emph{manifold hypothesis}.

In \cite[Theorem 2.13]{Dunne2023:SpreadDimension} we showed that for $(M, d_M)$ a submanifold $M\subset \R^k$ where $d_M$ is the Euclidean metric inherited from $\R^k$, if the limit $\lim\limits_{t\rightarrow \infty}\mathcal{G}\sigma_{d_M}(t)$ exists then
\begin{equation}\label{eq:1}
\mathcal{G}\sigma_{d_M}(t) \rightarrow \dim(M) \quad\text{ as }\quad t\rightarrow \infty
\end{equation}
where $\dim(M)$ is the topological dimension of $M$.

We then showed how the discrete and continuous spread dimension are related \cite[Theorem 3.9, Theorem 3.12]{Dunne2023:SpreadDimension}: we showed that for an increasing a sequence of finite subspaces $(X_n, d_n)\subset (M, d_M)$ satisfying $|X_n| \rightarrow \infty$ as $n \rightarrow \infty$, sampled in a sufficiently uniform way, then for each $t \in \Rplus$ we have
\begin{equation}\label{eq:2}
\mathcal{G}\sigma_{d_n}(t) \rightarrow \mathcal{G}\sigma_{d_M}(t) \quad\text{ as }\quad n \rightarrow\infty
\end{equation}

Taken together, (\ref{eq:1}) and (\ref{eq:2}) provide the theoretical basis for a new method of intrinsic dimension estimation, particularly when data is assumed to satisfy the manifold hypothesis. For large enough samples, an estimate for the intrinsic dimension of a dataset can be inferred from the graph of $\mathcal{G}\sigma_{d_n}(t)$ -- see Section \ref{sec:Conclusion}, or \cite[Section 4]{Dunne2023:SpreadDimension} for a more detailed discussion and examples with real and synthetic data.

The accuracy of $\mathcal{G}\sigma_{d_n}(t)$ as an intrinsic dimension estimator for a particular sample $X_n \subset M$ depends in part on the sample size $|X_n|$ being large enough, and hence the practical use of spread dimension as an intrinsic dimension estimator is closely tied to its computational feasibility.

\section{Approximating the Spread Dimension}

In this section we introduce the $S$-pseudo spread dimension of a metric space as an approximation of the spread dimension, and derive an expression for its standard error.

For each element $x\in X$ of a finite metric space $(X,d)$ we can define the function $\psi_{x}(t):\Rplus \rightarrow \Rplus$
\begin{equation*}\label{eq:PseudoSpread_x} \psi_{x}(t) = \frac{1}{\sum\limits_{y \in
X}e^{-td(x,y)}}|X|
\end{equation*}

For each $t\in \Rplus$, and each non-empty subset $S = \{x_1, x_2,...,x_k \} \subset X$, let $\Psi_S(t)$ denote the multiset
\begin{equation*}\label{eq:PseudoSpreadArray}
\Psi_S(t) = \{ \psi_{x_1}(t), \psi_{x_2}(t),..., \psi_{x_k}(t)\}
\end{equation*}

\begin{defn} Let $(X,d)$ be a finite metric space and let $S \subseteq X$
be a non-empty subset. The \emph{$S$-pseudo spread} $\overline{\Psi}_S(t):
\Rplus \rightarrow \Rplus$ is the function defined as the arithmetic mean of $\Psi_S(t)$, that is
\begin{equation*}\label{eq:PseudoSpread} \overline{\Psi}_{S}(t) =
\sum\limits_{x \in S}\frac{1}{\sum\limits_{y \in
X}e^{-td(x,y)}}\frac{|X|}{|S|}
\end{equation*} 
\end{defn}

Note that by construction, when $S=X$ the $S$-pseudo spread is equal to the spread, that is, $\sigma_d(t) = \overline{\Psi}_X(t)$.

\begin{defn}
For $(X,d)$ a finite metric space and $S \subseteq X$ a non-empty subset, the \emph{$S$-pseudo spread dimension} is the function $\mathcal{G}\overline{\Psi}_S(t): \Rplus \rightarrow \Rplus$, the instantaneous growth of the $S$-pseudo spread.
\end{defn}

Just as with the spread, when $S=X$ the $S$-pseudo spread dimension 
coincides with the spread dimension.

It will be useful to derive an exact 
expression for $\mathcal{G}\overline{\Psi}_S(t)$, and for this we will introduce the following quantities. For each $x \in X$, by direct calculation we have

\begin{equation}\label{eq:varphi_is_deriv} \frac{\dee}{\dee t} \psi_x(t) =
\frac{\sum\limits_{z \in
X}d(x,z) e^{-td(x,z)}}{\big(\sum\limits_{y \in
X}e^{-td(x,y)}\big)^2 }|X|
\end{equation}

Let $\varphi_{x}(t) = \frac{\dee}{\dee t} \psi_x(t)$, and for each non-empty subset $S = \{x_1, x_2,...,x_k \} \subseteq X$ let $\Phi_S(t)$ denote the multiset
\[
\Phi_S(t) = \{\varphi_{x_1}(t), \varphi_{x_2}(t),..., \varphi_{x_k}(t)\}
\]
and let $\overline{\Phi}_S(t)$ denote the arithmetic mean of $\Phi_S(t)$, that is
\[
\overline{\Phi}_S(t) = \sum\limits_{x \in S} \varphi_x(t)\frac{|X|}{|S|}
\]

We can now give a formula for $\mathcal{G}\overline{\Psi}_S(t)$ in terms of $\overline{\Phi}_S(t)$ and $\overline{\Psi}_S(t)$ which will be 
useful for deriving a formula for of the standard 
error in pseudo spread dimension as an approximation of the spread dimension.

\begin{thm}\label{thm:characterise_pseudo_spread_dimension}
For $(X,d)$ a finite metric space and $S \subseteq X$
a non-empty subset, the $S$-pseudo spread dimension has the 
following characterisaton
\begin{equation}\label{eq:pseudo_spread_dimension}
\mathcal{G}\overline{\Psi}_S(t) = t \frac{\overline{\Phi}_S(t)}{\overline{\Psi}_S(t)}
\end{equation}
\proof{Appliying the instantaneous growth formula 
(\ref{eq:definitionGrowthDimension}) to the definition of $S$-pseudo spread we have
\begin{align*}
\mathcal{G}\overline{\Psi}_S(t) &= \frac{t}{\overline{\Psi}_S(t)}\frac{\dee}{\dee t} \overline{\Psi}_S(t)
\end{align*}

Computing the derivative of $\overline{\Psi}_S(t)$ we have
\begin{align*}
\frac{\dee}{\dee t} \overline{\Psi}_S(t) &= \frac{\dee}{\dee t} \Big(\sum\limits_{x \in S} \psi_x(t) \frac{|X|}{|S|} \Big) \\ 
&= \sum\limits_{x \in S} \frac{\dee}{\dee t} \psi_x(t)\frac{|X|}{|S|}\\
&= \sum\limits_{x \in S} \varphi_x(t)\frac{|X|}{|S|}\\
&= \overline{\Phi}_S(t)
\end{align*}
as required. \QED}
\end{thm}

In \cite[Appendix A]{Dunne2023:SpreadDimension} we described algorithms to compute the spread and spread dimension of a finite metric space $(X,d)$ using array operations on the $|X|\times |X|$ array of values representing the distance matrix of the metric space. These algorithms have space complexity $\bigO(|X|^2)$.

With minor adaptations of those algorithms both $\overline{\Psi}_S(t)$ and $\overline{\Phi}_S(t)$ can be computed using array operations on the $|S|\times |X|$ array representing the partial distance matrix -- the values $d(s,x)$ for all $s \in S$ and $x \in X$.

By Theorem \ref{thm:characterise_pseudo_spread_dimension}, we see that to compute the $S$-pseudo spread dimension for each scale value $t \in \Rplus$ it is enough to compute $\overline{\Psi}_S(t)$ and $\overline{\Phi}_S(t)$, and therefore we can compute the $S$-pseudo spread dimension of $(X,d)$ in $\bigO(|S|\times |X|)$. Full implementation available at the author's GitHub\footnote{\url{https://github.com/dk-gh/spread_dimension}}

\subsection*{Deriving the Standard Error and Confidence Intervals}

For each $t \in \Rplus$, how well the $S$-pseudo spread approximates the spread can be conceptualised in terms of how well the \emph{sample mean} $\overline{\Psi}_S(t)$ approximates the \emph{population mean} $\overline{\Psi}_X(t)$, which can be expressed in terms of the variance of the values in $\Psi_S(t)$
\begin{equation*}
\var(\Psi_S(t)) =\frac{1}{|S|} \sum\limits_{x \in S}\big( \psi_{x} - \overline{\Psi}_S\big)^2 
\end{equation*}

Similarly, for each $t \in \Rplus$ we can calculate 
how well the sample mean $\overline{\Phi}_S(t)$ approximates the population 
mean $\overline{\Phi}_X(t)$ by computing the variance of $\Phi_S(t)$

\begin{equation*}
\var(\Phi_S(t)) =\frac{1}{|S|} \sum\limits_{x \in S}\big( \varphi_{x} - \overline{\Phi}_S\big)^2 
\end{equation*}

Using the characterisation of pseudo spread dimension from Theorem \ref{thm:characterise_pseudo_spread_dimension} we can compute an approximation of the standard error of the $S$-pseudo spread dimension using the following propagation of errors formula -- see, for example \cite[p. 36]{Lee2006:AnalyzingComplexSurvey} -- where for samples $a$ and $b$ with $A = \mean(a)$ and $B=\mean(b)$

\begin{equation}\label{eq:propagation_of_uncertainty}
\var\Big(\frac{A}{B}\Big) \approx \frac{A^2}{B^2} \Big(\, \frac{\var(a)}{A^2} + \frac{\var(b)}{B^2} - \frac{\covar(a,b)}{AB} \, \Big)
\end{equation}

Applying the propagation of errors formula (\ref{eq:propagation_of_uncertainty}) to the expression (\ref{eq:pseudo_spread_dimension}) we obtain the following approximation of the variance of the $S$-pseudo spread dimension for each value $t \in \Rplus$

\begin{equation*}\label{eq:propagation_of_uncertainty_applied}
\var\big(\mathcal{G}\overline{\Psi}_S(t)\big) \approx t^2\frac{\overline{\Phi}_S(t)^2}{\overline{\Psi}_S(t)^2} \Big(\, \frac{\var(\Phi_S(t))}{\overline{\Phi}_S(t)^2} + \frac{\var(\Psi_S(t))}{\overline{\Psi}_S(t)^2} - \frac{\covar(\Phi_S(t), \Psi_S(t))}{\overline{\Phi}_S(t)\overline{\Psi}_S(t)} \, \Big)
\end{equation*}

Where the covariance term is explicitly calculated as
\[
\covar(\Phi_S(t), \Psi_S(t)) = \frac{1}{|S|}\sum\limits_{x \in S}\big( \varphi_{x} - \overline{\Phi}_S\big)\big( \psi_{x} - \overline{\Psi}_S\big)
\]

The derived standard deviation is simply the square root of the variance $\SD =\sqrt{ \var\big(\mathcal{G}\overline{\Psi}_S(t)\big)}$. Following \cite[p. 422]{Freedman2007:Statistics} we take the standard deviation of the sample $\SD$ as an approximation of the standard deviation of the population giving an approximation of the standard error

\[
\SE \approx \frac{\sqrt{ \var\big(\mathcal{G}\overline{\Psi}_S(t)\big)}}{\sqrt{|S|}}
\]
and a $95\%$ confidence interval can be computed using the formula 
\begin{equation}\label{eq:CI95}
\mathcal{G}\overline{\Psi}_S(t) \pm 1.96 \times \SE
\end{equation}
where $1.96$ is an approximate value of the Z score for the $97.5\text{th}$ percentile point of the normal distribution.

Figure \ref{fig:swiss_10k} below shows an example computation using points sampled from the Swiss roll dataset\footnote{\url{https://scikit-learn.org/stable/modules/generated/sklearn.datasets.make_swiss_roll.html}}.

\begin{rem}
The same computation shown in Figure \ref{fig:swiss_10k} was repeated for 100 different randomly generated Swiss roll samples, and for each of these samples the pseudo spread dimension was computed across 200 separate scale values in the range $t \in [0, 15]$. The corresponding 95\% confidence intervals were computed for each $t$ using the formula (\ref{eq:CI95}). Across these 20,000 individual calculations, in $95.3\%$ of cases, the true pseudo spread dimension was within the calculated $95\%$ confidence interval, giving empirical validation of the derived standard error for pseudo spread dimension.
\end{rem}

\begin{figure}[H]
\centering
	\includegraphics[width=1\textwidth]{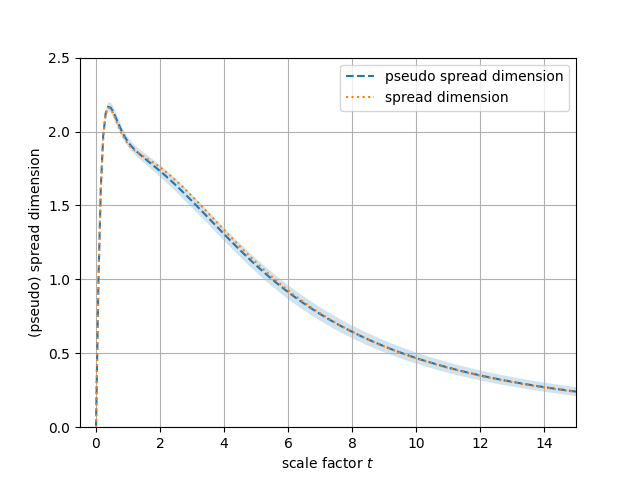}
	\caption{A Swiss roll dataset $X \subset \R^3$ with $|X| = 10,000$. Spread dimension is computed for a range of values $t \in [0,15]$, alongside the $S$-pseudo spread dimension computed for a random sample with $|S|=100$. Shaded area shows a 95\% confidence interval computed using the formula (\ref{eq:CI95}).}\label{fig:swiss_10k}
\end{figure}

\section{Conclusion}\label{sec:Conclusion}

Using spread dimension as an intrinsic dimension estimator requires computing the spread dimension across a range of scales and identifying a peak or plateau -- see \cite[Section 4]{Dunne2023:SpreadDimension} for full details. One important heuristic for applications as an intrinsic dimension estimator is that the longer a plateau around a particular value $n$, the stronger the indication that this is the true intrinsic dimension, i.e. that $\dim(X)=n$. 

\begin{figure}[H]
\centering
	\includegraphics[width=1\textwidth]{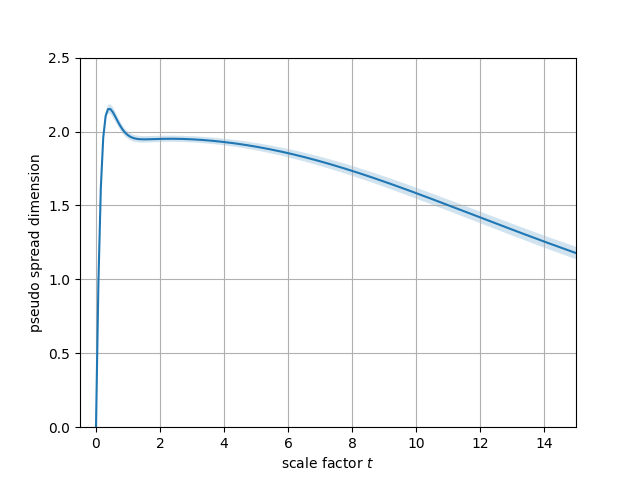}
	\caption{A Swiss roll dataset $X \subset \R^3$ with $|X| =100,000$. The $S$-pseudo spread dimension is computed for a sample with $|S|=100$ over a range of values $t \in [0,15]$. The shaded area depicts a 95\% confidence interval computed using the formula (\ref{eq:CI95}).}\label{fig:swiss_100k}
\end{figure}

For example, the Swiss roll example in Figure \ref{fig:swiss_10k} shows a peak at around $2$, whereas, in Figure \ref{fig:swiss_100k} we see the pseudo spread dimension computed for a larger sample of $100,000$ points gives a longer plateau near $2$, which corresponds with the intrinsic dimesion of the Swiss roll.

The introduction of pseudo spread dimension effectively makes computing the spread dimension of data more computationally feasible, which is an essential property for applying intrinsic dimension estimators in many real-world contexts involving large datasets.

Spread dimension has several other desirable properties as an intrinsic dimension estimation method:

\begin{itemize}
\item spread dimension does not require data to be of the form $X \subset \R^n$, as it is defined for a much broader class of metric spaces.
\item spread dimension is \emph{manifold adaptive} in the sense defined by Farahmand \cite{Farahmand2007:ManifoldAdaptive}, that is, the computational complexity of spread dimension of $X \subset \R^N$ does not depend on the dimension $N$ of the ambient space $\R^N$.
\item spread dimension can be applied to a dataset globally or locally -- see \cite[Section 4]{Dunne2023:SpreadDimension} for a detailed discussion.
\item spread dimension is inherently robust against multiscaling, as the intrinsic dimension is defined in terms of how the spread grows as the space is scaled -- see \cite[Section 4]{Dunne2023:SpreadDimension} for a detailed discussion.
\end{itemize}

\bibliography{bibliography}

\end{document}